%% file: chap-i.tex
\documentclass{amsart}

\include{headers}

\begin{document}
\title{Self Maps of $\P^1$ with Prescribed Ramification in
Characteristic $p$}
\author{Brian Osserman}
\begin{abstract}
Using limit linear series and a result controlling degeneration from
separable maps to inseparable maps, we give a formula for the number of
self-maps of $\P^1$ with ramification to order $e_i$ at general points
$P_i$, in the case that all $e_i$ are less than the characteristic. Unlike
the case of characteristic $0$, the answer is not given by Schubert 
calculus, nor is the number of maps always finite for distinct $P_i$, even 
in the tamely ramified case. However, finiteness for general $P_i$, obtained
by exploiting the relationship to branched covers, is a 
key part of the argument.
\end{abstract}
\thanks{This paper was partially supported by fellowships from the National
Science Foundation and Japan Society for the Promotion of Science.}
\maketitle

\section{Introduction}

We work throughout over an algebraically
closed field $k$ of characteristic $p$. For convenience, and because our 
results will be trivial in the case $p=2$, we assume $p>2$ throughout.
The question we wish to address in this paper is simply:

\begin{ques}\label{map-main-ques} Fix $n$ points $P_i$ on $\P^1$ and integers 
$e_i \geq 2$, with $\sum_i (e_i-1)=2d-2$, and $e_i \leq d$ for all $i$. 
How many separable self-maps of $\P^1$ of degree $d$ are there which ramify 
to order $e_i$ at the $P_i$, counted modulo automorphism of the image $\P^1$?
\end{ques}

The condition that $\sum_i
(e_i-1)=2d-2$ implies by the standard characteristic-$p$ Riemann-Hurwitz
formula that there are no solutions if any of the $e_i$ are not
prime to $p$, so we will assume throughout that all $e_i$ are prime to $p$
unless we specify otherwise.

In characteristic $0$, the number of maps is always finite, and when the
$P_i$ are general, this number is given combinatorially in terms
of Schubert calculus. In characteristic $p$, neither of these statement holds 
in general, and we explore the situation in a complete range of
characteristics, showing that the situation can be particularly pathological 
in low characteristics regardless of whether the ramification is tame or wild,
and ultimately solving the problem in mid range and higher characteristics.
Immediate motivation for this paper was given by applications to logarithmic 
connections with vanishing $p$-curvature on $\P^1$, and consequently 
Frobenius-unstable
vector bundles on curves of genus $2$, as well as the geometry of the
generalized Verschiebung; see \cite{os10}, \cite{os11}, and \cite{os9}.
However, the main question addressed is sufficiently fundamental that a
range of applications may be expected; in characteristic $0$,
generalizations of this question are applicable to computation of the 
Kodaira dimension of moduli spaces of curves (see \cite{lo1}) and the 
solution of the $A_N$ Bethe equation of XXX type (see \cite{m-v}).
Additionally, the results and techniques of this paper lead to a number of
new results on the existence and non-existence of branched covers of the
projective line; see \cite{os12}.

There is considerable
literature on our main question and its natural generalizations in
characteristic $0$, from
Eisenbud and Harris' original solution in the case of $\P^1$ in
\cite[Thm. 9.1]{e-h4}, to combinatorial formulas in the same cases by Goldberg
\cite{go1} and Scherbak \cite{sc1}, to formulas in the higher genus case of
Logan \cite[Thm. 3.1]{lo1} and the author \cite{os2}. However, the present 
work appears to be the first attempt to approach the problem for positive
characteristics.

Our basic technique is an adaptation of the limit linear series degeneration
argument of \cite{os2}, solving the problem first in the case of
three points, and then repeatedly letting ramification points come together 
to reduce inductively to this case. The main obstruction to carrying this
argument through is controlling potential degeneration of separable maps 
to inseparable maps. 

We now give some notation and terminology leading up to the statement of our
main theorem.

\begin{notn} When the answer to Question \ref{map-main-ques} is finite, we 
denote it by $N(\{(P_i, e_i)\}_i)$. We denote 
by $N_{\gen}(\{e_i\}_i)$ the value of $N(\{(P_i, e_i)\}_i)$ for general 
$P_i$.
\end{notn}

\begin{defn}
We distinguish three ranges of characteristic. We will refer to the {\bf
high characteristic} range to mean those characteristics for which $p>d$,
as well as characteristic $0$. The {\bf mid characteristic} range will be
characteristics for which $p \leq d$, but $e_i <p$ for all $i$. Finally, the
{\bf low characteristic} range will be characteristics for which $p \leq
e_i$ for some $i$. 
\end{defn}

We will see that high characteristics are uniformly
well-behaved with respect to our question, while low characteristics can be
extremely pathological, and the mid characteristics seem to be reasonably
well-behaved, but are considerably subtler than the high characteristics.

\begin{thm}\label{map-main}In the mid and high characteristics, we have the
following complete solution to our main question:

\begin{equation}N_{\gen}(\{e_i\}_i) = 
\!\!\!\!\!\!\!\!\!\!\!\!\!\!\!\!\!\!\!\!\!\!
\sum
_{\scriptsize \begin{matrix}d-e_{n-1}+1 \\ d-e_{n}+1\end{matrix} \leq d' \leq
\begin{matrix}d \\ p+d-e_{n-1}-e_n\end{matrix}} 
\!\!\!\!\!\!\!\!\!\!\!\!\!\!\!\!\!\!\!\!\!\!
N_{\gen}(\{e_i\}_{i \leq n-2}, e), \text{ with } e=2d'-2d+e_{n-1}+e_n-1
\end{equation}

\begin{equation}N_{\gen}(e_1, e_2, e_3)= \begin{cases}1 & p>d \\ 0 &
\text{otherwise} \end{cases}
\end{equation}
Further, for general points $P_i$ all of the relevant maps have no
non-trivial deformations.
\end{thm}

\begin{rem}We make a few observations: first, since the degree $d'$ is 
always no
greater than $d$, high characteristic will remain high under recursion.
Similarly, adding the two inequalities on the right, we find 
$e=2d'-2d+e_n+e_{n-1}-1 \leq d+(p+d-e_{n-1}-e_n)-2d+e_n
+e_{n-1}-1 = p-1$, so mid characteristic is also preserved (or becomes high)
under iteration.

Note that in the high characteristic range, we always have $p>d$, so the
answer becomes independent of characteristic: this is visibly true for the
second formula, and is true for the first formula because the inequality $d'
\leq p+d-e_{n-1}-e_n$ is subsumed by the inequality $d' \geq e =
2d'-2d+e_n+e_{n-1}-1$, or equivalently $d' \leq 2d-e_{n-1}-e_n+1$, which is
necessary for the number of maps $N_{\gen}(\{e_i\}_{i \leq n-2}, e)$ to be 
nonzero. Unsurprisingly, this
characteristic-independent formula is also the answer in characteristic $0$.
\end{rem}

We also remark that chronologically, the direct approach here
was not the first proof discovered of our formulas. That was obtained via
a correspondence with certain logarithmic connections on $\P^1$ together
with a theorem of Mochizuki, as outlined in 
\cite{os6}. The key step of the direct argument presented 
here, the analysis of
separable maps degenerating to inseparable ones, was derived via a careful
study of the corresponding situation with connections. 

We begin in Section 
\ref{s-map-translate} by translating the problem into a question on 
intersection of Schubert cycles in a Grassmannian. We exploit the 
relationship between ramified maps and branched covers in Section 
\ref{s-map-finite} to obtain some basic finiteness results including a 
ramified Brill-Noether-type theorem for $g^1_d$'s on $\P^1$ with specified 
ramification. We then apply this in Section \ref{s-map-3pt} to solve the 
base case of three ramification points. Section \ref{s-map-random} appears 
at first blush to be merely a couple of eccentric observations, including 
the pathology that when exactly one $e_i$ is greater than $p$, the number of 
maps can never be finite, but these 
observations play key roles in Section \ref{s-map-insep}, where we give a 
precise analysis of when a family of separable maps can degenerate to an 
inseparable map, and in Section \ref{s-map-degen}, where we finally prove 
our main theorem via the degeneration argument using limit linear series. 
Finally, in Section \ref{s-map-further} we explore some examples and 
further questions, and in Appendix \ref{s-map-moduli} we construct a scheme 
representing maps between a pair of fixed curves, with at least a certain 
specified ramification, but at points which are allowed to move; this 
scheme is the key idea in the proof of the Brill-Noether-type result of 
Section \ref{s-map-finite}, and is also used to generalize this result in
\cite{os3}.

The contents of this paper form a portion of the author's 2004 PhD thesis at 
MIT, under the direction of Johan de Jong.

\section*{Acknowledgements}

I would like to thank Johan de Jong for his tireless and invaluable
guidance. I would also like to thank Joe Harris, Dan Laksov, Astrid Giugni,
and Roya Beheshti for their helpful discussions.

\section{Translation to Schubert cycles}\label{s-map-translate}

In this section, we translate Question \ref{map-main-ques}
into a question on intersection of Schubert cycles on the projective 
Grassmannian $\G(1,d)$, and pin down some related notation. 
The translation is easy enough: a map (up to automorphism of the image) may
be represented explicitly by a $2$-dimensional space of polynomials of
degree $d$. It is then easy to see that a ramification condition of order
$e$ at a point $P$ corresponds to a Schubert cycle of codimension $e-1$,
which we denote by $\Sigma_{e-1}(P)$. 
Since we assumed $\sum_i (e_i-1) = 2d-2$, and our Grassmannian has dimension
$2d-2$, the expected dimension of
the intersection is therefore $0$. Pieri's formula will now give us the 
intersection product of our cycles,
yielding a hypothetical formula for the answer to our question. However,
there are several substantive issues to address. 

The first major issue is
whether or not the Schubert cycles will actually intersect transversely,
even for general choice of the $P_i$. Vakil \cite[Cor. 2.7 (a)]{va1} and Belkale \cite[Thm.
0.9]{be1} have recently independently shown that if the Schubert cycles are general, they will intersect transversely, but
it is not the case that general points on $\P^1$ will correspond 
to general Schubert cycles in $\G(1,d)$, so we cannot hope to
apply such general results. In fact, they correspond to osculating flags of the rational normal curve in $\P^d$. In characteristic $0$
properness of the intersection (that is, having
the expected dimension) for any choice of
distinct $P_i$ is straightforward, and we will reproduce the argument below
in order to analyze its implications in characteristic $p$. Transversality
for general choice of $P_i$ in characteristic $0$ is known, but more
involved (see \cite[Thm. 9.1]{e-h4}),
and means that Pieri's formula actually yields the correct number for
general choice of $P_i$. However, all of these statements fall apart
in characteristic $p$, as we will see shortly. 

The second issue to face is
that of base points: points of $\G(1,d)$ with base points correspond to
lower degree maps padded out by extra common factors in the defining 
polynomials. It is easy to see by Riemann-Hurwitz that base points cannot
occur for separable maps. In particular, in characteristic
$0$, or when $p>d$, the intersection of our Schubert cycles always actually
corresponds to the desired $g^1_d$'s. On the other hand, in general inseparable
maps can and will occur, frequently contributing an excess intersection. For
instance, in the case $d>p$, $e_i<p$, the Frobenius map will always
contribute a $\P^{d-p}$ to the intersection, with one point in $\G(1,d)$ for
every choice of a degree $d-p$ base point divisor.

These are the two issues which must be addressed in order to give an answer 
to the question. However, except in the base case of three points, we will 
not address them directly, as would be
required by an intersection-theoretic approach. We will rather take a 
different tack, looking at moduli of $g^1_d$'s with specified ramification
for certain degenerating families. Before continuing, we summarize:

\begin{prop} The answer to Question \ref{map-main-ques} is given by the
number of points corresponding to separable maps inside the
intersection $\cap _i \Sigma_{e_i-1} (P_i)\subset \G(1,d)$ of the Schubert 
cycles $\Sigma_{e_i-1}(P_i)$.
\end{prop}

\begin{warn}The equivalences between maps and linear series tend to become misleading in families; in
particular, if we have a linear series which develops base points in a
special fiber, there is no way to remove them globally to actually produce a
morphism. For our arguments, whenever we are working over a base other than a field, we will
always be dealing with linear series, even if we describe it as a ``family
of maps''. For the appropriate definitions (albeit in an overly generalized 
context), see \cite{os8}.
\end{warn}

\section{Finiteness Results}\label{s-map-finite}

We begin with a proposition whose argument is well-known in characteristic
$0$:

\begin{prop}In any characteristic, if $\sum _i (e_i -1) = 2d-2-c$ for some
$c \geq 0$, then any component of 
$\cap _i \Sigma_{e_i-1} (P_i)$ having dimension greater than $c$ must meet
the inseparable locus. In particular, in high characteristics, $\cap_i
\Sigma_{e_i-1} (P_i)$ always has the expected dimension $c$.
\end{prop}

\begin{proof} By Riemann-Hurwitz, if we had 
$\sum _i (e_i-1) > 2d-2$, then $\cap _i \Sigma_{e_i-1} (P_i)$ must consist 
entirely of inseparable maps. One can then induct on $c$, with a base
case of $c=-1$, and the induction step consisting of imposing a simple
ramification condition at a new point $P$. 
\end{proof}

The case $c=0$ is simply the full specification of a tame ramification 
divisor, so we restate:

\begin{cor}\label{map-highfinite}In high characteristics, there are finitely 
many self-maps of $\P^1$ with specified tame ramification divisor.
\end{cor}

The finite generation of fundamental groups of curves, together with some 
generalities on existence of moduli spaces of maps with certain 
ramification behavior, gives us a more substantive finiteness result 
than the previous proposition: 

\begin{thm}\label{map-gentame}Let $e_i$ be prime to $p$, and suppose $\sum_i (e_i-1)=2d-2$.
Then for a general choice of points $P_i$, we have that the set of maps from
$\P^1$ to $\P^1$ ramified to order $e_i$ at $P_i$, modulo automorphism of
the image: 
\begin{ilist}
\itm is finite;
\itm has no elements mapping any two of the $P_i$ to the same point.
\end{ilist}
\end{thm}

\begin{proof}By Theorem \ref{map-mr-main}, we have a moduli scheme
$\MR=\MR^d(\P^1, \P^1, \{e_i\}_i)$, with ramification and branching maps down 
to $(\P^1)^n$, and actions of $\Aut(\P_1)$ on both sides, with the action
on the domain being free. It is well-known
that given any specified tame branch locus, up to automorphism of the cover
there are only finitely many covers with the given degree and branching: 
this follows, for instance, from the finite generation of 
the tame fundamental group of
$\P^1$ minus the branch points.

Thus, each fiber of the branch morphism $\branch: \MR \rightarrow (\P^1)^n$ has 
only
finitely many $\Aut(\P^1)$ orbits, and is therefore of dimension at most 
$\dim \Aut(\P^1) = 3$. We conclude that the dimension of $\MR$ is at most $n+3$.
This immediately implies that a general fiber (in the sense of a fiber above
a general point of $(\P^1)^n$, making no hypotheses on dominance) of the 
ramification morphism $\ram: \MR \rightarrow (\P^1)^n$ can have dimension at 
most $3$. By the freeness of the $\Aut(\P^1)$-action on this side, this
completes the proof of (i). 

The proof of (ii) proceeds similarly: one sees that the locus $\MR'$ of 
maps in $\MR$ 
sending any two ramification points to the same branch point has dimension at 
most $n-1+3=n+2$, and the fibers of $\MR'$ under the ramification morphism are still $\Aut(\P^1)$-orbits, so we conclude that a general
fiber of the ramification morphism cannot contain any points of $\MR'$, completing the proof.
\end{proof}

\begin{rem} This finiteness theorem may be considered a first case in
positive characteristic of a Brill-Noether theorem with prescribed 
ramification, as in \cite[Thm. 4.5]{e-h1}. We show via deformation theory of
covers in \cite{os3} that one can generalize further in the $r=1$ case,
which gives in particular an instrinsically algebraic
and characteristic-$p$ proof of the previous theorem.
\end{rem}

\section{The Case of Three Points}\label{s-map-3pt}

While the general problem we wish to study becomes rather subtle in
characteristic $p$, the special case where we only have three ramification 
points is more tractable. This is fortuitous, since this
case will form the base case of our general induction argument. We begin 
by observing that in
this case, since each ramification index must be at most $d$, all three ramification points must map to distinct points. We can also show via elementary
observations that:

\begin{lem}\label{map-3ptPn}The intersection $\cap _i \Sigma_{e_i-1}(P_i)$ for
three points is isomorphic to $\P^m$ for some $m \geq 0$.
\end{lem}

\begin{proof} Set $P_1=0$, $P_2=\infty$, and $P_3 = 1$, and 
simultaneously fix bases $(F,G)$ (up to simultaneous scaling) for our 
$g^1_d$'s by requiring that 
$F$ vanish at $0$, $G$ vanish at $\infty$, and $F(1)=G(1)$. We may then 
verify the assertion directly by looking at the conditions imposed on the coefficients of $F$ and $G$ by the ramification conditions.
\end{proof}

We now show:

\begin{thm}\label{map-3pt}Let $P_1, P_2, P_3$ be three distinct points of $\P^1$, and $e_1,
e_2, e_3$ positive integers. Then we have: 
\begin{ilist}
\itm In any characteristic, $N(\{(P_i, e_i)\}_i)$ is finite, and is in 
fact always $0$ or $1$, being $0$ if and only if there is some inseparable 
$g^1_d$ of degree $d$ with the required ramification. Moreover, when 
$N(\{(P_i, e_i)\})=1$, the intersection is actually given 
scheme-theoretically by a single reduced point.
\itm Whenever $e_1$ and $e_2$ are less than $p$ and $d \geq p$, we have 
$N(\{(P_i, e_i)\}_i)=0$. Whenever $d<p$, we have $N(\{(P_i, e_i)\}_i)=1$.
\end{ilist}
\end{thm}

\begin{proof}We first deduce (ii) from (i): the second claim is trivial,
since if $d<p$, there can be no inseparable map of degree $d$. For the first
claim, because $e_1$ and $e_2$ are both less than $p$, any inseparable map
will satisfy the required ramification conditions at $P_1$ and $P_2$, and if 
we choose our map to be Frobenius, we can check directly that by adding base points at $P_3$ we can satisfy the last ramification condition as well.

For the proof of (i), we begin by noting that the intersection product in
question is always $1$: indeed, since all the Schubert cycles in question
are special, this follows immediately by applying 
Pieri's formula and then the
complementary-dimensional cycle intersection formula (see
\cite[duality theorem, p. 271]{fu1}. 
Next, the separable locus must be finite, since for three points on $\P^1$ there are no moduli, so we can apply Theorem \ref{map-gentame}.
Finally, by Lemma \ref{map-3ptPn}, our 
intersection is a $\P^m$, and is in particular connected. If 
it is $0$-dimensional, we are done, since we get a 
single reduced point which must clearly correspond to either a separable or 
inseparable map. On the other hand, if it is positive dimensional, by the closedness of the inseparable locus and the finiteness of the separable locus, we find that all the maps are inseparable. 
\end{proof}

To rephrase a slightly special case of the second part of the theorem, we 
have:

\begin{cor}Suppose we are in the situation of the preceding theorem, and
$e_1, e_2 <p$. Then a separable map of the specified ramification exists if 
and only if $d<p$.
\end{cor}

\begin{rem}This corollary certainly doesn't hold if we drop the hypothesis
that at least two ramification indices be less than $p$, as may be seen by considering the example of $x^n$.
\end{rem}

\section{Some Theorems and Pathologies}\label{s-map-random}

In this section, we make observations on what happens when some $e_i$ are
replaced by $p-e_i$, and we also find that when exactly one $e_i$ is greater
than $p$, there can never be a finite number of maps with the specified
ramification.

\begin{notn} Let $f$ be a separable map $f$ between smooth proper 
curves $C$ and $D$. Then the {\bf different} $\delta$ of $f$ is defined 
to be the divisor on $C$ associated to the 
skyscraper sheaf obtained as the cokernel of the natural map 
$f^* \Omega^1_D \hookrightarrow \Omega^1_C$.
\end{notn}

We have the following amusing and occasionally useful lemma.

\begin{lem}Fix $e_i$ all less than $p$, and points $P_i$ on $\P^1$. Given
an $f$ ramified to order $e_i$ at $P_i$, and with $f(P_1)\neq f(P_2)$, we
can associate an $\hat{f}$ ramified to order $p-e_i$ at $P_i$ for $i=1,2$
and $e_i$ at $P_i$ for $i>2$, and with $\hat{f}(P_1) \neq \hat{f}(P_2)$. 
This association is defined uniquely on
equivalence classes modulo automorphism of the image $\P^1$, and induces an
bijection on such equivalence classes.

In particular, if $e'_i$ are any integers obtained from the $e_i$ by 
repeatedly replacing pairs of indices $e_i, e_j$ with $p-e_i, p-e_j$ while 
holding the others fixed, we have
$N_{\gen}(\{e_i\}_i)=N_{\gen}(\{e'_i\}_i)$.
\end{lem}

\begin{proof}
For convenience, we assume that $f$ is unramified at infinity. 
By composing $f$ with an automorphism of the image $\P^1$,
we may write it (uniquely up to scalar) as 
$F/G=(x-P_1)^{e_1}F'/(x-P_2)^{e_2}G'$. If we multiply
through by $(x-P_2)^p/(x-P_1)^p$, we get the new function
$\hat{f}=(x-P_2)^{p-e_2}F'/(x-P_1)^{p-e_1}G'$. Since we obtained it from the 
old one
by multiplying by an inseparable function, one checks that the different is
unaffected away from $P_1$ and $P_2$. Since we assumed all $e_i<p$, it follows that the new function and old function have the same 
ramification away from $P_1$ and $P_2$, and {\it a priori} infinity. On the 
other
hand, it is clear that the ramification at $P_1$ and $P_2$ is now $p-e_1$
and $p-e_2$, and it is easy to check that the new degree of the function 
allows for no new ramification at infinity. This operation is visibly
invertible and well-defined up to automorphism equivalence, and in 
particular gives the desired bijection.

For the second statement, we just induct on pairs of $e_i, e_j$,
making use of the fact that by Theorem \ref{map-gentame}, for $P_i$ 
general, none of our maps for either of the two
relevant choices of ramification indices send any two of the
$P_i$ to the same point.
\end{proof}

The main usefulness of this rather eccentric fact is summarized in the 
following, to be applied later on:

\begin{cor}\label{map-reduction}To calculate $N_{\gen}(\{e_i\}_i)$ completely in the mid and high
characteristic range, it suffices to do so either when all but at most one
of the $e_i$ are less than $p/2$, or when all the
$e_i$ are odd. Moreover, it suffices to prove Theorem \ref{map-main} in only
either of these two cases.
\end{cor}

\begin{proof}The first statement follows trivially from the previous
corollary. The second is simply a matter of noting that for any given number
of points, the parity of the sum of the $e_i$ is determined by the integrality of
$d$. All the $e_i$ being odd always gives the correct 
parity, and if any $e_i$ are even, an even number of them must be.

To get the final assertion, we have to show that the formulas proposed in
Theorem \ref{map-main} are unaffected by replacing a pair $e_i$ and $e_j$ with 
$p-e_i$ and $p-e_j$. We will show this by induction, with $n=3$ as the base
case. For convenience, we repeat the formulas in question:

$$N_{\gen}(\{e_i\}_i) = 
\!\!\!\!\!\!\!\!\!\!\!\!\!\!\!\!\!\!\!\!\!\!
\sum
_{\scriptsize \begin{matrix}d-e_{n-1}+1 \\ d-e_{n}+1\end{matrix} \leq d' \leq
\begin{matrix}d \\ p+d-e_{n-1}-e_n\end{matrix}} 
\!\!\!\!\!\!\!\!\!\!\!\!\!\!\!\!\!\!\!\!\!\!
N_{\gen}(\{e_i\}_{i \leq n-2}, e), \text{ with } e=2d'-2d+e_{n-1}+e_n-1$$

$$N_{\gen}(e_1, e_2, e_3)= \begin{cases}1 & p>d \\ 0 &
\text{otherwise} \end{cases}$$

We begin with the three point case. Taking into account the additional inequalities $e_i \leq d$ for all $i$, if we substitute
$d=\frac{e_1+e_2+e_3-1}{2}$, we find our
inequalities may be rewritten as:
$$\begin{matrix}e_1-e_2+1 \\ e_2-e_1+1 \end{matrix} \leq e_3 \leq
\begin{matrix} e_1+e_2-1 \\ 2p-1-e_1-e_2 \end{matrix}.$$
Replacing any two $e_i$ by $p-e_i$ will simply permute these inequalities.

Next, since the proposed recursive equation is not visibly symmetric in the 
$e_i$, there are three cases to consider: first, $i=n-1$, $j=n$; second, 
$i,j<n-1$; and
finally, $i<n-1, j\geq n-1$. We first note that when we replace $e_i, e_j$ by $p-e_i, p-e_j$, 
the degree $d$ changes to $d+p-e_i-e_j$.

In the first case, one checks that under the substitutions, the inequalities for $d'$ are simply permuted, and $e$ remains unchanged.
For the second and third cases, we will want to write the inequalities for
$d'$ as equivalent inequalities for $e$. We find:
$$\begin{matrix}e_n-e_{n-1}+1 \\ e_{n-1}-e_n+1 \end{matrix} \leq e \leq
\begin{matrix}e_n+e_{n-1}-1 \\ 2p-1-e_{n-1}-e_n\end{matrix}$$
In particular, these inequalities depend only on $e_n$ and $e_{n-1}$, and
not on $d$. Hence, in the second case $e$ ranges through the same values after substituting for $e_i$ and $e_j$, so we need only use the induction
hypothesis to conclude the desired result.

Finally, in the third case we assume for convenience that $i<n-1$ and $j=n$.
In this case, when we substitute into our inequalities for $e$, we get
$$\begin{matrix}p-e_n-e_{n-1}+1 \\ e_{n-1}+e_n+1 -p\end{matrix} \leq e \leq
\begin{matrix}p-e_n+e_{n-1}-1 \\ p-1-e_{n-1}+e_n\end{matrix}$$
which then gives us that $p-e$ satisfies precisely the same inequalities
that $e$ did originally. Thus, each term in the new recursive formula
corresponds to a unique one in the old one by replacing $e$ with $p-e$ and
$e_i$ by $p-e_i$, so once again the induction
hypothesis gives us the desired result.
\end{proof}

We end with a rather surprising observation illustrating that even tame
ramification can have very pathological behavior in low characteristics:

\begin{prop}\label{map-path}Suppose $e_1>p$ but still prime to $p$, and $e_i<p$ for all $i>1$. 
Then if one map exists with ramification $e_i$ at $P_i$, infinitely many do.
In particular, if the $P_i$ are general, no maps exist with ramification
$e_i$ at $P_i$.
\end{prop}

\begin{proof}Without loss of generality, we may assume that $P_1$ is the
point at infinity, and that our function maps infinity to infinity, so that
it is given by $F/G$, with $\deg F-\deg G = e_1$. Now consider the family 
of functions $F/G - tx^p$, where $t \in k$. Because $e_1>p$, the
ramification at infinity is unaffected. On the other hand, since $x^p$ is
regular away from infinity and inseparable, the different is
unchanged on the affine part, and since we assumed that all $e_i<p$ for
$i>1$, we find that the ramification is unaffected everywhere, giving us an
infinite family of maps, clearly not related by automorphism, all with the
same ramification. For $P_i$ general we know that there can be at most
finitely many maps with specified tame ramification by Theorem
\ref{map-gentame}, so we conclude that there cannot be any such maps at all.
\end{proof}

The following corollary will not be used later, but seems worth mentioning:

\begin{cor}\label{map-path-wild} Suppose that the $P_i$ are
general, and we have a map $f$ with $e_i <p$ for all $i>1$, but 
$e_1 = mp$ wild with $m > 1$. Then the different of $f$ at $P_1$ is greater
than $2(m-1)p$. 
\end{cor}

\begin{proof} If we again put $P_1$ at infinity and write $f=F/G$,
subtracting off some multiple of $x^{mp}$ will force the degree to drop,
and leave all $e_i$ for $i>1$ unchanged. The index $e_1$ may not drop 
(this can only
happen if the degree of $F$ drops at least $mp$ below the degree of $G$), but 
if it remains wild we can iterate, and since the degree drops each time,
$e_1$ must eventually become tame. By our proposition, this new tame index,
which we denote by $e'_1$, would have to be less than $p$. If we denote
the new degree by $d'$, one sees that $d-d' > (m-1)p$, considering separately the cases that the degree of $F$ remained greater than $G$ or dropped below that of $G$. The corollary then follows from Riemann-Hurwitz.
\end{proof}

\begin{ex}To demonstrate that the statement of Proposition \ref{map-path} is not
vacuous, we note it is not difficult to write down a concrete example.
Indeed, the family of functions $x^{p+2}+tx^p+x$ is easily seen to satisfy our requirements.
\end{ex}

\begin{rem} In standard examples of linear series existing only for special configurations in characteristic $0$, the space of linear series with the prescribed
ramification is supported over a maximal-dimensional subspace of
$\cM_{g,n}$, or equivalently, a general configuration where such a linear
series exists has only finitely many. In particular, in none of the standard
examples is the expected dimension non-negative. It is not clear 
whether or not this must always be the case in characteristic $0$, but 
here we have an example where this fails to hold in characteristic $p$.
\end{rem}

\section{Specialization to Inseparable Maps}\label{s-map-insep}

The ultimate goal will be to solve the map-counting problem in mid and high
characteristics by repeatedly
letting points come together. The main obstacle to this is understanding
when a family of separable maps can have an inseparable map as its limit.
We provide an answer to this question which may seem 
unmotivated, and indeed arose from careful examination of the situation in
the very different and at first glance totally unrelated setting of
certain connections with vanishing $p$-curvature, as discussed in 
\cite{os6}.

Our main result is:

\begin{thm}\label{map-insep}Let $A$ be a DVR containing its residue field $k$
and with uniformizer $t$, and $f_t$ be a family of maps of degree $d$ from 
$\P^1$ to $\P^1$ over $\Spec A$ (more precisely, a linear series on $\P^1_A$)
whose generic fiber is tamely ramified along sections $P_i$ with all $e_i<p$, 
and whose special
fiber is inseparable. We further assume that the $P_i$ stay away from
infinity. Then if the limit of the $P_i$ in the special fiber is denoted by
$\bar{P}_i$, we have:
\begin{ilist}
\itm If the $\bar{P}_i$ are distinct, they are in a special 
configuration allowing the existence of separable maps of degree
$d+m-1$
ramified to order $e_i$ at $\bar{P}_i$, and $2 m - 1$ at infinity;
\itm If $\bar{P}_j =\bar{P}_{j'}$ with $e_j+e_{j'}<p$, and the other
$\bar{P}_i$ distinct, then the $\bar{P}_i$ are in a special
configuration allowing separable maps of degree $d+m-1-b+\epsilon$, 
ramified to order $e_i$ at the $\bar{P}_i$ for $i\neq j,j'$, 
to $e_j+e_{j'}-2b-1$ at $\bar{P}_j=\bar{P}_{j'}$, and to
$2 m -1$ at infinity; 
\end{ilist}
in either case, $m$ is some integer with $p \leq m \leq d$, and in
the second case $b$ is a non-negative integer less than $(e_j+e_{j'}-1)/2$.
\end{thm}

\begin{proof}The main idea of the proof is not dissimilar to the basic
operation of applying fractional linear transformations to be able to factor
out a power of the uniformizer if one is given a family of maps degenerating
to a constant map. However, in this case we will apply a fractional linear
transformation with inseparable coefficients; this will behave similarly,
but will not preserve the degree of the map, and also does not appear to
work readily in nearly the generality of the constant case.

We work for the most part explicitly with pairs of polynomials and their
differents, only dealing with common factors at the end to translate to 
rational functions and ramification indices.
We can write $f_t$ as $F/G$, with $F, G \in A[x]$, and
have no common factors. We denote by $F_0$ and $G_0$ the polynomials
obtained from $F$ and $G$ by setting $t=0$, and by $\bar{F}_0$ and
$\bar{G}_0$ the inseparable polynomials obtained by canceling the common
factors of $F_0$ and $G_0$; since $F,G$ represent a linear series of
dimension $1$, we may further assume that they were chosen so that $F_0,
G_0$ defines a non-constant function. Then let $H_1$ and $H_2$ be inseparable
polynomials of degree strictly less than $\bar{F}_0$ and $\bar{G}_0$
respectively, such that $\bar{F}_0 H_2 - \bar{G}_0 H_1 = 1$ (this is
possible by dividing the exponents $\bar{F}_0$ and $\bar{G}_0$ by $p$, 
applying Euclid's algorithm in $k[x]$, and multiplying all exponents by 
$p$). We now
construct a new family $\tilde{F}/\tilde{G}$ over $\Spec A$ as follows: if we
denote by $\nu$ the map from $A[x]$ to itself which simply factors out
common powers of $t$, then $\tilde{F} := \nu(F \bar{G}_0 - G \bar{F}_0)$,
and $\tilde{G} := F H_2 - G H_1$. It is easy to check that applying an
inseparable fractional linear transformation to $F/G$ will change
$(dF)G-F(dG)$ by the determinant of the transformation, so it follows that
$(d\tilde{F})\tilde{G}-\tilde{F}(d\tilde{G})$ is the same as $(dF)G-F(dG)$, 
but with a positive power of $t$ factored out. 

At $t=0$,
we note that  since we had $\bar{F}_0 H_2 - \bar{G}_0 H_1 = 1$, $\tilde{G}$ 
is made up precisely of the common factors of $F_0$ and $G_0$, of which there 
can be at most $d- \deg f_0 \leq d-p$. Since we removed a positive power of 
$t$ from $(dF)G-F(dG)$, if we still have an inseparable limit, we can repeat 
the process as many times as necessary to remove all the powers of $t$ and
obtain a separable limit. Each time we do, the degree of the denominator at 
$t=0$ clearly remains at most $d-p$. We thus end up with a family
$\tilde{F}/\tilde{G}$ which over the generic fiber has the same different
as $F/G$ away from infinity. If we let $K$ be the fraction field of $A$, we
also note that we must have that the ideal generated by
$\tilde{F},\tilde{G}$ in $K[x]$ is the same as that generated by $F,G$.
Since $F,G$ had no common factors over $K$, it follows that 
$\tilde{F}, \tilde{G}$ have no common factors either. Now, since we have no
common factors, we find by considering differents that away from $t=0$
(that is, at the generic fiber), 
$\tilde{F}/\tilde{G}$ has the same ramification as $F/G$ except possibly at
infinity, since all the $e_i$ were specified to be less than $p$. 

Denote by $\tilde{F}_0,\tilde{G}_0$ the polynomials obtained from
$\tilde{F},\tilde{G}$ at $t=0$. We claim that 
the different of $(\tilde{F}_0,\tilde{G}_0)$ away from 
infinity is the limit of the different of $f_t$, and in particular has
degree $2d-2$.
Indeed, this follows from our hypothesis that the $P_i$ stay away from
infinity, because when the limit is separable, the limit of the different 
is the different of the limit, with orders adding when points come together.
Next, from our hypotheses, the different in the limit is less than $p$
except at infinity, and we need not worry
about wild ramification.  Denote by $e_{\infty}$ the ramification index of 
$\tilde{F}_0/\tilde{G}_0$ at infinity, and suppose that $p|e_{\infty}$. In
this situation, we can replace $\tilde{F}_0$ by subtracting off an
appropriate multiple of $x^{e_{\infty}}\tilde{G}_0$, which will decrease 
$e_{\infty}$
without affecting the ramification away from infinity. Repeating as
necessary, we can require that $\tilde{F}_0/\tilde{G}_0$ be tamely ramified
at infinity (hence everywhere), without changing its behavior away from
infinity, or the degree of $\tilde{G}_0$.

If we denote the greater of the degrees of $\tilde{F}_0, \tilde{G}_0$
by $\tilde{d}$, we then have that
$$2\tilde{d}-2=2d-2+e_\infty-1.$$ 
We find in particular that
$\tilde{d}\geq d$, so since the degree of $\tilde{G}_0$ was strictly
less than $d$, it follows that $\tilde{d}$ is simply the degree of 
$\tilde{F}_0$, and $e_\infty = \tilde{d} - d_0$, where $d_0$ is the degree of 
$\tilde{G}_0$. Write $d_0=d-m$ for some $m$, where from our earlier bound on
the degree of our denominators, we know that $m \geq p$. We then find that 
$\tilde{d}=d+m-1$, and $e_{\infty}=2m-1>p$.

Finally, we translate back into the language of maps, by removing common 
factors. Because vanishing indices cannot drop under
specialization, and the different of the limit is the limit of the
different, it then follows that at any point $P$, we can acquire base points
only if more than one ramification section converges to $P$, which is to say
only if $P=\bar{P}_j=\bar{P}_{j'}$. In this situation, we let $b$ be the 
number of common factors of $\tilde{F}_0, \tilde{G}_0$ at 
$\bar{P}_j = \bar{P}_{j'}$. It now follows immediately that removing these
common factors, we have constructed a function as claimed in the statement
of the theorem.
\end{proof}

Putting the theorem together with Proposition \ref{map-path}, and noting
that there are only finitely many possibilities for $m$ and $b$, we 
conclude:

\begin{cor}\label{map-no-insep} In the situation of the preceding theorem,
if the $\bar{P}_i$ are general, there cannot be any $f_t$ as described,
having an inseparable limit. 
\end{cor}

\section{The Degeneration Argument}\label{s-map-degen}

We complete the proof of Theorem \ref{map-main} in this section via a
degeneration argument. The basic situation we will consider is the family
arising as follows:

\begin{figure}
\centering
\input{reduced.pstex_t}
$$\text{An explicit family of smooth $\P^1$'s degenerating to a node, with
$n$ sections.}$$
\end{figure}

Inside of $\P^1 \times \P^1$, take the family of hyperbolas 
given in affine coordinates by $xy=t$, degenerating at $t=0$ to the union 
of the $x$-axis and $y$-axis. For each $t \neq 0$, we get a smooth $\P^1$, 
and fix isomorphisms between them by projecting to the $y$-axis.
Choose an isomorphism between our abstract $\P^1$ and the $y$-axis sending $P$
to the node; we can now speak of $P_1, \dots, P_{n-2}$ as well as $P$ as 
fixed points on the $y$-axis and simultaneously on all the $\P^1$'s in our 
family; they are (constant) sections of our family. Now, choose any two 
points $P_{n-1}^0$ and $P_n^0$ on the $x$-axis away from $0$, and define 
sections $P_{n-1}^t$ and $P_n^t$ similarly via projection from our family 
to the $x$-axis rather than the $y$-axis. Under our fixed trivialization 
of the smooth fibers of the family, these sections both tend towards 
the section defined by $P$ (see figure). We will consider this as a family
$X$ over $\Spec k[t]$, and write $X_t$ for the associated local
family over $\Spec k[t]_{(t)}$.

We briefly review the main concepts of the theory of limit linear series 
as it relates to our situation. See \cite{os8} for general
definitions and, where applicable, proofs. On any non-singular fiber of our 
family, we
know that a map to $\P^1$ (modulo automorphism of the image) corresponds to a 
$g^1_d$ on that fiber; we see that given a $g^1_d$ on the family away from
the special fiber, we can obtain a $g^1_d$ on either the $x$- or $y$-axis
simply by projecting all fibers to the appropriate choice of axis. This pair
gives the associated Eisenbud-Harris limit series on the nodal fiber; 
we have vanishing sequences $a^x_i$ and $a^y_i$ for $i=0,1$ at the node, 
and the degree of the induced map on the $x$-axis (resp., $y$-axis) 
is at most $d-a^x_0$ (resp., $d-a^y_0$), with the ramification index 
of the map at the node given by $a^x_1 - a^x_0$ (resp., $a^y_1 - 
a^y_0$). We have the inequalities $a^x_0 + a^y_1 \geq d$, $a^x_1 + a^y_0 
\geq d$; the data of a pair of $g^1_d$'s on the components with vanishing 
sequences satisfying these inequalities is in fact the definition of an
Eisenbud-Harris limit series, and we say that a given limit series is refined 
if these are both equalities. 

It is easy to see that for arbitrary $a^x_i$ 
and $a^y_i$ satisfying the necessary inequality, the space of $g^1_d$'s on
the $x$-axis will have total required ramification at least $2(d-a^x_0)-2$,
and similarly for the $y$-axis, so in particular by Riemann-Hurwitz if the 
limits are
separable, we immediately conclude that they must form a refined limit 
series, and they cannot have any additional base points, so the corresponding 
maps must have degrees precisely $d-a^x_0$ and $d-a^y_0$, with ramification 
index $a^x_1-a^x_0 = a^y_1-a^y_0$ at the node.

Given these observations, our general theory, and specifically 
\cite[Thm. 5.3]{os8}, gives us:

\begin{thm}Associated to our families $X$ and $X_t$, and any 
choice of ramification indices $e_i$ such that $\sum _i (e_i-1)=2d-2$, 
are schemes $G^1_d:= G^1_d(X, \{(P_i, e_i)\}_i)$ and $\tilde{G}^1_d := 
G^1_d(X_t, \{(P_i, e_i)\}_i)$, with the latter obtained from the former by
base change, and the fibers parametrizing (limit) linear series with the
required ramification on the fibers of $X$ and $X_t$. We also have
open subschemes $G^{1,\sep}_d$ and $\tilde{G}^{1,\sep}_d$ parametrizing limit
series which are separable when restricted to every component of every
fiber. Over $t=0$, $G^{1,\sep}_d$ (equivalently, $\tilde{G}^{1,\sep}_d$)
parametrizes simply Eisenbud-Harris limit series, and contains only refined
limit series.
\end{thm}

\begin{proof} Most of this is immediate from \cite[Thm. 5.3]{os8}. The fact 
that $G^{1,\sep}_d$ parametrizes Eisenbud-Harris series on the special fiber
follows from \cite[Cor. 6.8]{os8} together with the assertion that the 
only separable Eisenbud-Harris limit series are refined, which we observed
above. 
\end{proof}

Given this language, we can readily apply \cite[Cor. 6.12]{os8} to obtain:

\begin{cor}\label{map-finiteflat}With the notation of the above theorem, if 
$P_1, \dots, P_{n-2}$ and $P$ are chosen 
generally, and $e_{n-1}+e_n<p$, then $\tilde{G}^{1,\sep}_d$ 
is finite etale over $\Spec k[t]_{(t)}$. In particular, it has the same
number of points, all reduced, in the geometric generic and special fibers,
and the fibers of $G^{1,\sep}_d$ have the same number of points at general 
$t$ as at $t=0$.
\end{cor}

\begin{proof} First, the assertion on the fibers of $G^{1,\sep}_d$ for 
general $t$ follows immediately from the statement on
$\tilde{G}^{1,\sep}_d$, together with the fact that $\tilde{G}^{1,\sep}_d$ 
is obtained from $G^{1,\sep}_d$ simply by localization of the base around 
$t=0$.

Next, to obtain the desired statement on $\tilde{G}^{1,\sep}_d$, we need 
only verify that the conditions (I)-(III) of \cite[Cor. 6.12]{os8} are 
satisfied: first, that every separable
Eisenbud-Harris limit series on the special fiber is refined; second, that
the scheme of separable Eisenbud-Harris limit series on the special fiber
consists of a finite number of reduced points; and third, that if $A$ is a
DVR, any $A$-valued point of $\tilde{G}^1_d$ mapping flatly to 
$\Spec k[t]_{(t)}$
and being separable at the generic point is also separable on the closed
point. Condition (I) is satisfied even without
the generality hypothesis, as stated in the above theorem. 

Condition (III) is
for the most part simply an application of Corollary \ref{map-no-insep}; 
indeed, given an $A$-valued point of $\tilde{G}^1_d$ flat over $\Spec
k[t]_{(t)}$, 
projection to the $y$-axis would give a family of $g^1_d$'s on $\P^1$ with
ramification sections specializing to the $P_1, \dots, P_{n-2}, P$, which are
general by hypothesis. Then Corollary \ref{map-no-insep} says that if the 
family is generically separable, it must remain separable on the special 
fiber. It
remains to see that the same holds if we project to the $x$-axis. 
For this, considering the different we note that the vanishing sequence on 
the $y$-axis at the node will satisfy $a^y_0+a^y_1-1=e_{n-1}+e_n-2$, and in 
particular $a^y_1 < p$. On the other hand, $a^x_0+a^y_1 \geq d$, so since
$a^x_0$ is the number of base points acquired on the $x$-axis, the degree on
the $x$-axis is less than or equal to $d-a^x_0 \leq a^y_1 <p$, and 
we also cannot have an inseparable limit along the $x$-axis, giving
condition (III).

Lastly, we prove the validity of condition (II) by induction on $n$. 
The basic observation is because the space of refined Eisenbud-Harris limit
series may be viewed simply as a disjoint union over all vanishing sequences
satisfying $a^x_i+ a^y_{1-i} =d$ of the products of the schemes
parametrizing $g^1_d$'s with appropriate ramification on each component, it
suffices to see that these latter are made up of reduced points. It is easy
to see that as the vanishing sequences vary, if we simply remove the base
points $a^x_0$ and $a^y_0$, we will have the same ramification index $e$ at
the node on each component, the degrees on each component will be such that
the expected dimension (taking $e$ into account as well as the $e_i$) will
be zero, and $e$ will vary arbitrarily given this constraint, together
with the constraint that the degrees on each component be at most $d$. In
particular, it suffices to see that for points chosen generally, the 
scheme of separable $g^1_d$'s in the $(n-1)$-point and $3$-point cases always 
consist of
a finite number of reduced points, and by induction on the statement of our 
corollary, it is enough to see this in the $3$-point case, which we have
conveniently already handled in Theorem \ref{map-3pt}. 
\end{proof}

We are now ready for:

\begin{proof}[Proof of Theorem \ref{map-main}]
First, we may assume that $e_{n-1}+e_n<p$, thanks to Corollary
\ref{map-reduction}.
By Corollary \ref{map-finiteflat}, for all our points
chosen generally, and a general choice $t$, $G^{1,\sep}_d$ has the same 
number of points over that particular $t$ as it does over $t=0$. This sets
up a simple recursion formula to calculate $N_{\gen}(\{e_i\}_i)$: the number
will be given by the number over the special fiber, which is the sum over
all choices $e$ of ramification index at the node of
$N_{\gen}(e,e_{n-1},e_n) N_{\gen}(\{e_i\}_{i<n-1},e)$.

We recall
that the formula we wanted to prove for Theorem \ref{map-main} (the second
formula having already been handled by Theorem \ref{map-3pt}) was
$$N_{\gen}(\{e_i\}_i) = 
\!\!\!\!\!\!\!\!\!\!\!\!\!\!\!\!\!\!\!\!\!\!
\sum
_{\scriptsize \begin{matrix}d-e_{n-1}+1 \\ d-e_{n}+1\end{matrix} \leq d' \leq
\begin{matrix}d \\ p+d-e_{n-1}-e_n\end{matrix}} 
\!\!\!\!\!\!\!\!\!\!\!\!\!\!\!\!\!\!\!\!\!\!
N_{\gen}(\{e_i\}_{i \leq n-2}, e), \text{ with } e=2d'-2d+e_{n-1}+e_n-1$$
and that in the proof of Corollary \ref{map-reduction} we saw that the above
inequalities for $d'$ were equivalent to the following inequalities on $e$:
$$\begin{matrix}e_n-e_{n-1}+1 \\ e_{n-1}-e_n+1 \end{matrix} \leq e \leq
\begin{matrix}e_n+e_{n-1}-1 \\ 2p-1-e_{n-1}-e_n\end{matrix}$$

We first show that the above inequalities for $e$ are precisely the
range for which $N_{\gen}(e, e_{n-1}, e_n)=1$. But with Theorem \ref{map-3pt} at
our disposal, this is a trivial observation, since $e_n-e_{n-1}+1 \leq e$,
$e_{n-1}-e_n+1 \leq e$ and $e \leq e_n + e_{n-1}-1$ are precisely the
inequalities insuring that the ramification indices are less than the degree
of the map, and $e \leq 2p-1-e_{n-1}-e_n$ insures that the degree is less
than $p$. Finally, we need to know that the degree on the three-point
component will be less than $d$. This degree will be given by
$\frac{e+e_{n-1}+e_n-1}{2}$, so we find that {\it a priori}, we need $e \leq
2d-e_{n-1}-e_n+1$. However, we note that the right hand side is actually
$2d'-e$, so this inequality is equivalent to $e \leq d'$, and we needn't
include it with the conditions, as if it is violated we will have
$N_{\gen}(\{e_i\}_{i \leq n-2}, e)=0$, and there will be no contribution to
the sum. This completes the proof of our main theorem.
\end{proof}

\section{Examples and Further Questions}\label{s-map-further}

We first apply our main theorem in the case of four points. For a given $d'$
as in Theorem \ref{map-main}, $N_{\gen}(e_1, e_2, e)=1$ if $e_1, e_2, 
e \leq d'$ and 
$p > d'$, and $N_{\gen}(e_1, e_2, e)=0$ otherwise. Rewriting this condition
in terms of $d'$, we get the bounds $e_1 \leq d'$, $e_2 \leq d'$, 
$d' \leq 2d - e_3 - e_4 +1$, $d' \leq
p-1$, and including these bounds for $d'$ along with those of Theorem
\ref{map-main}, simply by substracting the various bounds
for possible values of $d'$ we obtain:

\begin{cor}\label{map-4pts} The number $N_{\gen}(\{e_i\}_i)$ of self-maps 
of $\P^1$ of degree $d$
in characteristic $p$, ramified to orders $e_1, \dots, e_4$ at four general
points, with each $e_i < p$ and $2d-2 = \sum _i (e_i -1)$, and counted modulo 
automorphism of the image, is given by the formula
$$N_{\gen}(\{e_i\}_i) = \min\{\{e_i\}_i, \{d+1-e_i\}_i, \{p- e_i\}_i,
\{p-d-1+e_i\}_i\},$$
or equivalently,
$$N_{\gen}(\{e_i\}_i) = \min\{\{e_i\}_i, \{d+1-e_i\}_i\}-\max\{0, d+1-p\}.$$
Further, all of these maps are without any nontrivial deformations.
\end{cor}

\begin{ex} We explore an example which demonstrates all the basic behaviors
we have described so far, and may be solved explicitly: maps of degree 3,
with four simple ramification points. We may assume without loss of
generality that $P_1=0, P_2 = \infty, P_3 = 1,$ and we let $P_4$ be a
general parameter $\lambda$. We see immediately that our four ramification
points must have distinct images, so we may further specify that our maps fix 
$P_1$ and $P_2$, from which we deduce that they are of the form
$f=\frac{x^2(ax+b)}{x+c}$, with $a,b,c$ all nonzero. Since we did not specify
that $P_3$ be fixed, we have one remaining degree of freedom, and may set
$b=1$. Now, if we consider the zeroes and poles of $df$, we can calculate
directly that our possible maps satisfy $2c= 2a \lambda$ and $1+3ac =
-(1+\lambda)2a$, which in characteristic $\neq 2$ means $c$ is determined by
$a$ and $\lambda$, and $a$ satisfies $3 \lambda a^2 +2(1+\lambda)a+1=0$. In
characteristic $3$, we get a unique (separable) solution, while in
characteristics $0$ or $p>3$, we get two solutions for general $\lambda$. We
find that these solutions come together when $1-\lambda +\lambda^2=0$. 
Finally, in characteristic 3, we also see that the unique solution $f =
\frac{x^3+(1 + \lambda)x^2}{(1+\lambda)x+ \lambda}$ specializes 
to an inseparable one when $\lambda$ goes to $-1$.
\end{ex}

\begin{rem}\label{map-path-rem} In the context of covers with prescribed branching, there is a very general result on lifting to characteristic $0$; see \cite[11, Proof of Prop. 5.1]{b-l-r2}. One can use this and the three-transitivity of $\Aut(\P^1)$ to prove a similar result for lifting $g^1_d$'s on $\P^1$ with prescribed ramification at three points, and one
might be tempted to conjecture that one can generalize to
arbitrary numbers of points and tame ramification on $\P^1$. However, Proposition \ref{map-path} shows that such a statement cannot hold even in this case, since
by virtue of Corollary \ref{map-highfinite} only finitely many of the 
infinitely many constructed maps would be able to
lift to characteristic $0$. However, it may still be true that one can 
prove such a statement in 
the high and mid characteristic ranges. Indeed, we are able to prove this in
more generality in \cite[Cor. 5.5]{os8}, subject to an
expected-dimension hypothesis. In particular, in our case of self-maps of 
$\P^1$, we can conclude thanks to Theorem \ref{map-gentame} that lifting to 
characteristic $0$ is always
possible for tame ramification indices and general ramification points, or
in high characteristic and arbitrary ramification points, and we have 
reduced the mid-characteristic case down to the question of finiteness for
arbitrary distinct ramification points. 
\end{rem}

We conclude with some further questions. We could reasonably start with
remaining questions about the case of $\P^1$, including:

\begin{ques}Is it true that for a given $d$ and $e_i$, the number of maps is
either always finite or always infinite as the $P_i$ are allowed to move? 
Can we prove that it is always finite in the mid-characteristic case?
\end{ques}

\begin{ques}What happens in low characteristic when more than one
ramification index is greater than $p$? 
\end{ques}

\begin{ques}What can we say about the dimension of spaces of wildly ramified 
maps? When do wildly ramified maps exist for general ramification points?
\end{ques}
 
This first question is partially answered by applying results of Mochizuki in
\cite{os6}, and similar arguments may be expected to give a complete affirmative answer. The last question is explored further in \cite{os3}. We briefly examine some examples of the second question. The following lemma is trivial.

\begin{lem} Let $f,g$ be rational functions on a smooth curve $C$, and $P \in C$ a point where both $f$ and $g$ are defined. Then $f+g$ has ramification index at $P$ at least as large as $f$ if and only if the ramification index of $g$ at $P$ is at least as large as that of $f$. 

In particular, in the case that $C=\A^1$ and $f$ is a tamely ramified polynomial of degree $d$, with $d>p$ and prime to $p$, we have that $f+g$ is another polynomial of degree $d$, with the same ramification as $f$, if and only if $g$ is an inseparable polynomial of degree less than $d$, and with ramification index greater than that of $f$ at every point of $\A^1$.
\end{lem}

\begin{ex} We conclude that if $e_1=d$ is strictly between $2p$ and $3p$, for some choices of $e_i$ we have a positive-dimensional family of maps with the given ramification, while for others we only have finitely many. For instance, if we have $p<e_2<2p$, but $e_i<p$ for all $i>2$, we obtain a unique inseparable polynomial (up to fractional linear transformation) of degree $2p$ ramified to order $2p$ at $P_2$, which gives us a one-dimensional family of maps with the desired ramification. But if $e2>2p$, or if $e_3$ is also between $p$ and $2p$, no inseparable polynomial of the specified form exists, so our separable map is necessarily unique. Furthermore, these examples are non-vacuous: examples of the first and second may be obtained simply by $x^d-x^{e_2}$, for $e_1-e_2 \neq p$, while $2x^{2p+1}-x^{2p}-2x^{p+1}$ is an example of the third.
\end{ex}

Lastly, one could ask the same questions about maps from
higher-genus curves to $\P^1$. These have been answered in the case of
characteristic $0$ in \cite{os2}, and the argument there would also apply in
characteristic $p$ given an appropriate generalization of Theorem
\ref{map-insep} to control the possibility of separable maps specializing to
inseparable maps. The case of higher-dimensional linear series is still
open as well, even in characteristic $0$. 

\appendix

\section{Moduli Schemes of Ramified Maps}\label{s-map-moduli}

The goal of this appendix is to construct moduli schemes of maps of curves
required to have at least given ramification, but at unspecified points.
Before we begin, we recall the well-known corollary of Grothendieck's work
on the Hilbert scheme:

\begin{thm}\label{map-mr-mor}Given $X$ and $Y$ two smooth, projective, geometrically connected
curves over a locally Noetherian scheme $S$ and a positive integer $d$, 
then the functor 
$\cMor^d_S(X,Y)$ parametrizing degree $d$ morphisms from $X$ to $Y$ over $S$
is representable
by a quasi-projective scheme. In particular, $\cAut_S(X)= \cMor^1_S(X,X)$ is
representable.
\end{thm}

\begin{proof}
Without the degree hypothesis, the functor is constructed in \cite[p.
221-20]{gr1} (where it is called $\Hom$) as an open subscheme of the 
Hilbert scheme via the graph
associated to a morphism. Now, if $\L$ and $\M$ are ample line bundles on
$X$ and $Y$, and $f$ a morphism of degree $d$, one checks that the Hilbert polynomial of the graph under
the projective imbedding of $X \times _S Y$ induced by $\pi_1^* \L \otimes
\pi_2 ^* \M$ is uniquely determined by $d$, so the $\Mor$ scheme is naturally a 
disjoint union over all $d$ of schemes representing $\Mor^d$, each of which 
is quasi-projective.

In order to see that $\cAut_S(X)=\cMor^1_S(X,X)$, we first note that for any
$d>0$, $\cMor_S^d(X,Y)$ consists entirely of scheme-theoretically surjective
morphisms. Indeed, given $f \in \cMor^d$, one checks by the criterion on flatness and fibers (see \cite[Thm. 11.3.10]{ega43}) that $f$ is faithfully flat, which implies scheme-theoretic surjectivity. Now, to see that
$\cAut_S(X) = \cMor^1_S(X,X)$, it suffices to note that in our situation,
one can check whether $f$ is a closed immersion on each fiber $f_s$ 
(see \cite[Prop. 4.6.7]{ega31}), so the desired assertion follows from
the well-known case of smooth curves over $S = \Spec k$ (see, for instance,
\cite[Cor. 4.4.9]{ega31}).
\end{proof}

We also have:

\begin{prop}With the notation of the preceding theorem, there exists an
open subscheme $\Mor^{d,\sep}_S(X,Y)$ of $\Mor^d_S(X,Y)$ parametrizing 
morphisms which are separable on every fiber.
\end{prop}

\begin{proof}Let $M:= \Mor^d_S(X,Y)$, $X_M$ and $Y_M$ be the pullbacks of $X$
and $Y$ to $M$, $\tilde{f}:X_M \rightarrow Y_M$ be the universal morphism of 
degree $d$, defined over $M$; we get an induced map $\tilde{f}^*
\Omega^1_{Y_M/M} \rightarrow \Omega^1_{X_M/M}$ of line bundles on $X_M$,
with the kernel giving the locus on $X_M$ where $\tilde{f}$ is ramified. The
complement is an open set, and its image in $M$ is clearly the locus of
separable maps; since $X$ is flat and of finite type over $S$, $X_M$ is flat
and of finite type over $M$, and in particular open, so we have constructed
an open subscheme of $M$ corresponding to separable maps, as desired.
\end{proof}

We also recall a standard construction involving the jet bundle, or bundle
of principal parts $\sP^n_{X/S}$, associated to an $S$-scheme $X$. The 
terminology and notation is not standard, however. 

\begin{defn}We define the {\bf $n$th cotangent bundle} $\Upsilon^n_{X/S}$
to be the kernel of the natural map $\sP^n_{X/S} \rightarrow \O_X$;
explicitly, consider $\O_X \otimes _{\O_S} \O_X$ as an $\O_X$-module via
left multiplication, and consider the natural map to $\O_X$ sending $a
\otimes b$ to $ab$. Then if we denote the kernel of this map by $\I_{X/S}$,
with the induced $\O_X$-module structure, $\Upsilon^n_{X/S} :=
\I_{X/S}/\I_{X/S}^{n+1}$.
\end{defn}

We recall:

\begin{prop}\label{map-mr-ega}With notation as in the preceding definition,
\begin{ilist}
\itm $\Upsilon^n_{X/S}$ is compatible with base change.
\itm On affine opens, $\I_{X/S}$ is generated by elements of the form
$a \otimes 1 - 1 \otimes a$, for $a \in \O_X$.
\itm If $X$ is smooth over $S$, $\Upsilon^n_{X/S}$ is locally free. 
\end{ilist}
\end{prop}

\begin{proof}
Compatibility with base change for $\sP^n_{X/S}$ is \cite[Prop.
16.4.5]{ega44}; because $\Upsilon^n_{X/S}$ is the kernel of a map
(clearly compatible with base change) to $\O_X$, and $\O_X$ is free, it
follows that $\Upsilon^n_{X/S}$ is compatible with base change.
(ii) is \cite[Lem. 0.20.4.4]{ega41}. Finally,
(iii) follows from the same statement for $\sP^n_{X/S}$, which is 
\cite[Prop. 17.12.4]{ega44}, since $\Upsilon^n_{X/S}$ is the
kernel of a surjective map from $\sP^n_{X/S}$ to $\O_X$ (in fact, this is
somewhat gratuitous, since the argument for $\sP^n_{X/S}$ works without
modification for $\Upsilon^n_{X/S}$).
\end{proof}

We now specify in full detail the functor we wish to represent. 

\begin{defn} Suppose we are given a pair of smooth, projective, 
geometrically connected curves $X,Y$ over a locally Noetherian base $S$, as
well as $n$ integers $e_i$, and $d \geq 1$. 
Then the functor $\cMR_S^d(X,Y,\{e_i\}_i)$ associates to any 
scheme $T$ over $S$
the set of separable morphisms $f$ from $X_T$ to 
$Y_T$ over $T$ of degree $d$, together with a choice of $n$ disjoint 
$T$-valued points 
$P_i$ of $X_T$, such that the fiber of $f(P_i)$ contains an $e_i$th-order
thickening of $P_i$ inside of $X_T$ for each $i$.
\end{defn}

Conceptually, this functor is the functor of maps $f$ of degree $d$ between 
$X$ and $Y$, together with points $P_i$ on $X$ which are (at least)
$e_i$th-order ramification points of $f$.

Our main result is:

\begin{thm}\label{map-mr-main}The functor $\cMR=\cMR^d_S(X,Y,\{e_i\}_i)$ is
representable by a scheme $\MR$. We also have 
the natural data of morphisms 
$\ram: \MR \rightarrow X^n$ and $\branch: \MR \rightarrow Y^n$ and 
actions of the group schemes $\Aut(X)$ and $\Aut(Y)$ on $\MR$ over $Y^n$ and 
$X^n$ respectively. Furthermore, $\Aut(Y)$ acts freely on $\MR$. 
\end{thm}

\begin{proof}First we note that all the assertions other than
representability can be verified simply on the functor level: the morphism
$\ram$ is the forgetful transformation which takes a point of $\MR$ and
remembers only the $P_i$; similarly, the morphism $\branch$ remembers the
$f(P_i)$, which are sections of $Y_T$. A point $g$ of $\Aut(Y)$ act on points 
of $\MR$ by sending $f$ to $g \circ f$ and leaving the $P_i$ fixed, and 
similarly $g \in \Aut(X)$ acts on $\MR$ by sending $f$ to $f \circ g$ and 
the $P_i$ to $g^{-1} P_i$, which fixes $f(P_i)$. The freeness of the
$\Aut(Y)$ action follows easily from the statement that any point of
$\MR$ corresponds to a scheme-theoretically surjective map, noted in the
proof of Theorem \ref{map-mr-main}. 

Clearly, we have a forgetful map from $\MR$ to
$M=\Mor^{d,\sep}(X,Y)$; since the latter is representable, it will enough 
to show that the map of functors is also representable. In fact, if we use 
the convention that $X_{M}^n$ denotes the product of $n$ copies of $X_{M}$ 
over $M$, the sections in the definition of our functor will allow us to 
describe $\MR$ as a 
closed subscheme of $X_{M}^n$ with the pairwise diagonals removed. We 
claim that it is 
enough to handle the case $n=1$: suppose we have done this case, and for
each $i$ let $\MR_i$ be the resulting scheme; then if we imbed the product of the $\MR_i$ as a closed subscheme of $X_{M}^n$, and remove pairwise diagonals, we get the desired scheme.

Since $X$ and $Y$ are smooth over $S$ by hypothesis,
$\Upsilon^{e-1}_{X/S}$ and $\Upsilon^{e-1}_{Y/S}$ are locally free,
so the kernel of any morphism $f^* \Upsilon^{e-1}_{Y/S} \rightarrow
\Upsilon^{e-1}_{X/S}$ is representable by a closed subscheme of $X$ over $S$,
and the following lemma completes the proof of our theorem.
\end{proof}

\begin{lem}Let $f:X \rightarrow Y$ be a morphism of separated $S$-schemes. 
Then there
is a natural map $f^* \Upsilon^{e-1}_{Y/S} \rightarrow \Upsilon^{e-1}_{X/S}$ 
such 
that for any $T$ over $S$, and any section $\sigma:T \rightarrow X_T$ we have: 

$(f^* \Upsilon^{e-1}_{Y/S} \rightarrow \Upsilon^{e-1}_{X/S})_{\sigma(T)} = 0$ if and
only if the fiber of $f_T$ over $f_T(\sigma(T))$ contains an $e$th order
thickening of $\sigma(T)$ inside $X_T$.
\end{lem}

\begin{proof}The map from $f^* \Upsilon^{e-1}_{Y/S} \rightarrow
\Upsilon^{e-1}_{X/S}$ is simply the one induced by $f^* \otimes f^*: 
f^{-1} \O_Y \otimes_{\O_S} f^{-1} \O_Y \rightarrow \O_X \otimes_{\O_S} \O_X$. 
Our assertion is local, so we immediately reduce to affines, and consider the
situation that $X_T = \Spec A$, $Y_T = \Spec B$, and $T = \Spec R$. Since
$X$ and $Y$ are separated over $S$, a section is a closed immersion, so we 
also denote by $I_\sigma$ the ideal corresponding to $\sigma(T)$ in $X_T$, 
and set $I_\sigma':= I_\sigma \otimes _R A \subset A \otimes _R A$. Now, 
the fiber of $f(\sigma(T))$ is given by 
$\Spec (A/I_\sigma \otimes _B A)$, cut out in $X_T$ by the ideal of $A$ 
generated by $(f_T^* B) \cap I_\sigma$. The fiber contains an 
$e$th order thickening of $\sigma_T$ if and only if this ideal is contained 
in $I_\sigma^e$, which is to say, if and only if $(f^*_T B)\cap I_\sigma
\subset I_\sigma^e$. We claim that this is equivalent to $1 \otimes ((f^*_T
B) \cap I_\sigma) \subset (I_\sigma', A \otimes I_\sigma^e)$: indeed, this
follows immediately from the fact that $(A \otimes _R A) / (I_\sigma', A
\otimes I_\sigma^e) \cong A/I_\sigma \otimes R A/I_\sigma^e \cong 
A/I_\sigma^e$. 

On the other hand,
$(f^* \Upsilon^{e-1}_{Y/S} \rightarrow \Upsilon^{e-1}_{X/S})_{\sigma(T)} 
= 0$ if and only if $(f_T^* \Upsilon^{e-1}_{Y_T/T} \rightarrow
\Upsilon^{e-1}_{X_T/T})_{\sigma(T)} = 0$, by Proposition \ref{map-mr-ega} (i),
and this is equivalent to the assertion that $f_T^* \I_{Y_T/T}$ is contained
in the ideal generated by $I'_\sigma$ and $\I_{X_T/T}^e$. Thus, the proof of
the lemma is reduced to the following two assertions: first, that
modulo $I'_\sigma$, we have $f_T^* \I_{Y_T/T} = 1 \otimes
((f^*_T B) \cap I_\sigma)$; and second, that  
$(I'_\sigma, \I_{X_T/T}^e) = (I'_\sigma, A \otimes _R I_\sigma^e)$. By 
Proposition \ref{map-mr-ega} (ii), $\I_{X_T/T}$ and $\I_{Y_T/T}$ are 
generated by elements of the form $a \otimes 1 - 1 \otimes a$ and 
$b \otimes 1 - 1 \otimes b$, respectively. The main observation is that
because $I_\sigma$ is the ideal of a section, to generate these ideals it
suffices to restrict to $a$ with $a \in I_\sigma$, and to $b$ with 
$f_T^* b \in I_\sigma$. The latter immediately gives the
first assertion. The second assertion also follows easily, since we then
have that $\I_{X_T/T}^e$ is generated by products $\prod_{j \leq e} (i_j
\otimes 1 - 1 \otimes i_j)$, which are equivalent modulo $I'_{\sigma}$ to 
elements of $1 \otimes I_\sigma^e$, as desired.
\end{proof}

\begin{rem}The $\Aut(X)$ action is not free in
general, often having a non-trivial finite sub-group scheme stabilizing any
given morphism. However, it is easy enough to see that the stabilizer of any
$k$-valued point $f \in \MR$ is in fact a finite group scheme. Indeed, in 
this case, we may as well set $S= \Spec(k)$. Since $\Aut(X)$ is a 
finite-type group 
scheme, the stabilizer will likewise be a finite-type group scheme 
over $k$, and it thus suffices to show that it consists of only finitely 
many $\bar{k}$-valued points. 
Now, an automorphism of $X_{\bar{k}}$ is determined on the
generic point, and will have to fix $K(Y_{\bar{k}})$ inside 
$K(X_{\bar{k}})$ in order to
fix $f$; since $K(Y_{\bar{k}})$ is a finite subfield of $K(X_{\bar{k}})$, 
the the relevant automorphism group is finite, and we conclude the desired assertion. 
\end{rem} 

\bibliographystyle{hamsplain}
\bibliography{hgen}
\end{document}

%% file: headers.tex
\usepackage{amssymb,euscript,amsmath, mathrsfs}
\usepackage[dvips]{graphicx}
\usepackage[dvips]{color}

\newcounter{ENUM}
\newcommand{\itm}{\item}
\newenvironment{ilist}{\renewcommand{\theENUM}{\roman{ENUM}}\renewcommand{\itm}{\addtocounter{ENUM}{1}\item[(\theENUM)]}\begin{itemize}\setcounter{ENUM}{0}}{\end{itemize}}

\newcommand{\margh}[1]{}

\input xy
\xyoption{all}
\CompileMatrices

\def\A{{\mathbb A}}
\def\P{{\mathbb P}}

\def\G{{\mathbb G}}

\def\L{{\mathscr L}}
\def\M{{\mathscr M}}

\def\O{{\mathscr O}}
\def\I{{\mathscr I}}

\def\sP{{\mathscr P}}

\def\cM{{\mathcal M}}

\def\cMor{\mathcal{M}or}
\def\cAut{\mathcal{A}ut}

\def\cMR{\mathcal{MR}}

\def\gen{\operatorname{gen}}

\def\Mor{\operatorname{Mor}}
\def\Hom{\operatorname{Hom}}

\def\Aut{\operatorname{Aut}}
\def\sep{\operatorname{sep}}
\def\ram{\operatorname{ram}}
\def\branch{\operatorname{branch}}

\def\Spec{\operatorname{Spec}}

\def\MR{\operatorname{MR}}

\numberwithin{equation}{section}
\newtheorem{thm}{Theorem}[section]
\newtheorem{prop}[thm]{Proposition}
\newtheorem{lem}[thm]{Lemma}
\newtheorem{cor}[thm]{Corollary}

\theoremstyle{definition}
\newtheorem{defn}[thm]{Definition}
\newtheorem{ques}[thm]{Question}
\newtheorem{ex}[thm]{Example}

\theoremstyle{remark}
\newtheorem{notn}[thm]{Notation}
\newtheorem{rem}[thm]{Remark}
\newtheorem{warn}[thm]{Warning}

%% file: reduced.pstex_t
\begin{picture}(0,0)%
\includegraphics{reduced.pstex}%
\end{picture}%
\setlength{\unitlength}{1973sp}%
\begingroup\makeatletter\ifx\SetFigFont\undefined%
\gdef\SetFigFont#1#2#3#4#5{%
  \reset@font\fontsize{#1}{#2pt}%
  \fontfamily{#3}\fontseries{#4}\fontshape{#5}%
  \selectfont}%
\fi\endgroup%
\begin{picture}(7224,7224)(2389,-7573)
\put(6151,-7261){\makebox(0,0)[lb]{\smash{{\SetFigFont{10}{12.0}{\familydefault}{\mddefault}{\updefault}{\color[rgb]{0,0,0}$P_1$}%
}}}}
\put(6151,-5086){\makebox(0,0)[lb]{\smash{{\SetFigFont{10}{12.0}{\familydefault}{\mddefault}{\updefault}{\color[rgb]{0,0,0}$P_{n-2}$}%
}}}}
\put(8476,-4336){\makebox(0,0)[lb]{\smash{{\SetFigFont{10}{12.0}{\familydefault}{\mddefault}{\updefault}{\color[rgb]{0,0,0}$P_n$}%
}}}}
\put(7201,-4336){\makebox(0,0)[lb]{\smash{{\SetFigFont{10}{12.0}{\familydefault}{\mddefault}{\updefault}{\color[rgb]{0,0,0}$P_{n-1}$}%
}}}}
\put(6226,-4336){\makebox(0,0)[lb]{\smash{{\SetFigFont{10}{12.0}{\familydefault}{\mddefault}{\updefault}{\color[rgb]{0,0,0}$P$}%
}}}}
\put(6226,-6061){\makebox(0,0)[lb]{\smash{{\SetFigFont{10}{12.0}{\familydefault}{\mddefault}{\updefault}{\color[rgb]{0,0,0}$\vdots$}%
}}}}
\end{picture}%